\documentclass[12pt,a4paper]{article}
\usepackage{amsmath,amssymb,amsfonts}
\usepackage[latin1]{inputenc}

\def\Bbb{\mathbb}

\title{\bf On the variance of the digits of $1/p$}

\author{Kurt Girstmair\\
Institut f\"ur Mathematik \\
Universit\"at Innsbruck   \\
Technikerstr. 13/7        \\
A-6020 Innsbruck, Austria \\
Kurt.Girstmair@uibk.ac.at}

\date{}

\makeatletter
\let\@@maketitle=\maketitle
\def\maketitle{\def\thispagestyle##1{\relax}\@@maketitle}
\makeatother
\newtheorem{theorem}{Theorem}

\newtheorem{corollary}{Corollary}

\def\BE{\begin{equation}}
\def\EE{\end{equation}}
\def\BD{\begin{displaymath}}
\def\ED{\end{displaymath}}
\def\BA{\begin{array}}
\def\EA{\end{array}}
\def\BEA{\begin{eqnarray*}}
\def\EEA{\end{eqnarray*}}
\def\BI{\bibitem}

\def\Z{\Bbb Z}
\def\Q{\Bbb Q}

\def\phi{\varphi}
\def\EPS{\varepsilon}

\def\MB{\mbox}
\def\LD{\ldots}
\def\OV{\overline}

\def\DIV{\,|\,}
\def\NDIV{\, \nmid \,}

\def\BQ{``}
\def\EQ{'' }

\def\MN{\medskip\noindent}
\def\STOP{\hfill$\Box$}

\def\BPS{B_{\psi}}
\def\BPSB{B_{\OV{\psi}}}
\def\BCHPS{B_{\chi\psi}}
\def\RE{\rm{Re}}

\newcommand{\btop}[2]{\genfrac{}{}{0pt}{1}{#1}{#2}}

\def\LS#1#2{ \left( \frac{#1}{#2} \right) }

\begin{document}
\maketitle
\normalsize

\begin{abstract}

Let $p>3$ be a prime and $b\ge 2$ an integer such that $p$ does not divide $b$. Then $1/p$ has a periodic digit expansion with respect to the basis $b$. The length $q$ of the period
is the (multiplicative) order of $b$ mod $p$. In the case $q=p-1$ a formula for the variance of the digits of a period was given previously. This formula involves a Dedekind sum.
We determine the variance in the case $q=(p-1)/2$.
If $p\equiv 3$ mod 4 a Dedekind sum and the class number of $\Q(\sqrt{-p})$ occur in the respective formula. If $p\equiv 1$ mod 4, the formula may be much more complex since
it involves linear combinations of (possibly many) products of two Bernoulli numbers attached to odd characters.

\end{abstract}

%%%%%%%%%%%%%%%%%%%%%%%%%%%%%%%%%%%%%%%%%%%%%
\section*{1. Introduction and results}
%%%%%%%%%%%%%%%%%%%%%%%%%%%%%%%%%%%%%%%%%%%%%

Let $p\ge 3$ be a prime, $b\ge 2$ an integer such that $p\NDIV b$.
Then
\BE
\label{1.1}
  \frac 1p=\sum_{j=1}^{\infty}c_jb^{-j},
\EE
where the numbers $c_j\in\{0,1,\LD,b-1\}$ are the {\em digits} of $1/p$ with respect to the basis $b$. It is well-known that the sequence of the digits is periodic and
that $(c_1,\LD,c_q)$ is a period, $q$ being the (multiplicative) order of $b$ mod $p$; see \cite{Gi1}.

In the paper \cite{Gi2} we determined the {\em variance}
\BD
 \sigma^2=\frac 1q\sum_{j=1}^q (c_j-m)^2
\ED
of the period $(c_1,\LD,c_q)$, where $m$ is the {\em mean value}
\BE
\label{1.1.1}
  m=\frac 1q\sum_{j=1}^q c_j,
\EE
in the case $q=p-1$.
Indeed, we obtained
\BE
\label{1.2}
\sigma^2=\frac{2bs(p,b)}{p-1}+\frac{(b-1)(bp-3b+p+3)}{12(p-1)},
\EE
where
\BD
s(p,b)=\sum_{k=1}^{b-1}((k/b))((pk/b))
\ED
is the classical Dedekind sum; for its definition see \cite[formula (1)]{RaGr}.

In this paper we treat the case $q=(p-1)/2$. Here Dedekind sums together with additional arithmetical data appear. The case $p\equiv 3$ mod 4 is, as a rule, much simpler than the case $p\equiv 1$ mod 4.
In what follows let $h(d)$ denote the class number of the field $\Q(\sqrt d)$. Moreover, let $T(b,p)$ denote the right-hand side of (\ref{1.2}).

\begin{theorem} % Theorem 1 %%%%%%%%%%%%%%%%%%%%%%%%%%%%%%%%%%%%%%%%%%%%%%%%%%%%%%%%%%%%%%%%%%%%%%%%
\label{t1}
Let $p\equiv 3$ mod $4$, $p>3$, and suppose that $b$ has the order $q=(p-1)/2$ mod $p$. Then
\BD
   \sigma^2=T(b,p)-\frac{(b-1)^2}{(p-1)^2}\,h(-p)^2.
\ED
\end{theorem} %%%%%%%%%%%%%%%%%%%%%%%%%%%%%%%%%%%%%%%%%%%%%%%%%%%%%%%%%%%%%%%%%%%%%%%%%%%%%%%%%%%%%%

In other words, the variance in the case of Theorem \ref{t1} differs from (\ref{1.2}) by the summand $-(b-1)^2h(-p)^2/(p-1)^2$.

\MN
{\em Example.} Let $p=151$, $b=10$. Here $s(p,10)=3/5$ and $h(-p)=7$. Theorem \ref{t1} gives $\sigma^2=5046/625=8.0736$.

The case $p\equiv 1$ mod $4$ requires additional notation. Let $\chi=\LS -p$ denote the Legendre symbol, which is an even Dirichlet character mod $p$. We also need odd Dirichlet characters $\psi$ mod $b$
and products $\chi\psi$ defined by $\chi\psi(k)=\chi(k)\psi(k)$. The characters $\chi\psi$ are odd Dirichlet characters mod $pb$. By $X_b^-$ we denote the set of odd characters mod $b$.
For $\psi\in X_b^-$, let $\BPS$ and $\BCHPS$ be the Bernoulli numbers
\BD
  \BPS=\frac 1b\sum_{k=1}^{b-1}\psi(k)k\:\MB{ and}\: \BCHPS=\frac 1{pb}\sum_{k=1}^{pb-1}\chi\psi(k)k.
\ED
The complex conjugate character of $\psi$ is denoted by $\OV{\psi}$. By $\phi$ we denote, as usual, Euler's function. Recall that $T(b,p)$ is the right-hand side of (\ref{1.2}).

\begin{theorem} % Theorem 2 %%%%%%%%%%%%%%%%%%%%%%%%%%%%%%%%%%%%%%%%%%%%%%%%%%%%%%%%%%%%%%%%%%%%%%%%%%%%%%%%%%%%
\label{t2}
If $p\equiv 1$ mod $4$ and $b$ has the order $q=(p-1)/2$ mod $p$, then
\BE
\label{1.3}
   \sigma^2=T(b,p)+\frac{2b}{p-1}\sum_{d\DIV b}\frac{\chi(d)}{\phi(d)}\sum_{\psi\in X_d^-}\OV{\psi}(p)\BCHPS\BPSB.
\EE
\end{theorem} %%%%%%%%%%%%%%%%%%%%%%%%%%%%%%%%%%%%%%%%%%%%%%%%%%%%%%%%%%%%%%%%%%%%%%%%%%%%%%%%%%%%%%%%%%%%%%%%%%

If $\psi\in X_d^-$, $d\DIV b$, is a primitive character, the $\BCHPS$ and $\BPSB$ do not vanish. If $\psi$ is imprimitive, let $\psi'$ be the primitive character mod $f$, $f\DIV d$, that induces $\psi$.
Then
\BE
\label{1.4}
 \BPS=\prod_{l\DIV d}(1-\psi'(l))B_{\psi'} \MB{ and } \BCHPS=\prod_{l\DIV d}(1-\chi\psi'(l))B_{\chi\psi'},
\EE
where $l$ runs through the prime divisors of $d$; see \cite[p. 274]{Sz}. In this case the Bernoulli numbers $\BCHPS$ and $\BPSB$ may vanish.

\MN
{\em Examples.} 1. Let $p\equiv 1$ mod $4$ and $b=10$. There are two primitive odd characters $\psi_5$, $\OV{\psi_5}$
mod 5, where $\psi_5$ is defined by $\psi_5(2)=i$, $i=\sqrt{-1}$. These two characters induce two imprimitive characters mod 10.
In addition, let $\EPS\in\{\pm 1\}$. By means of the formulas (\ref{1.4}), we obtain the following corollary.

\begin{corollary} % Corollary 1 %%%%%%%%%%%%%%%%%%%%%%%%%%%%%%%%%%%%%%%%%%%%%%%%%%%%%%%%%%%%%
\label{c1}
In the above setting,
\BD
  \sigma^2=T(10,p)+\begin{cases} \EPS\frac 6{p-1}\RE(B_{\chi\psi_5}(-3+i)) & \MB{ if }p\equiv \EPS\MB{ mod } 10;\\
                          \EPS\frac 2{p-1}\RE(B_{\chi\psi_5}(7+i)) & \MB{ if }p\equiv 3\EPS\MB{ mod } 10.
                   \end{cases}
\ED
\end{corollary} %%%%%%%%%%%%%%%%%%%%%%%%%%%%%%%%%%%%%%%%%%%%%%%%%%%%%%%%%%%%%%%%%%%%%%%%%%%%%

For instance, in the case $p=373$ we have $T(10,p)=1019/124$. Since $B_{\chi\psi_5}=10-4i$, the second case of the corollary (with $\EPS=1$) gives
$\sigma^2=1019/124+37/93=3205/372$.

2. Let $p\equiv 1$ mod 4 and $b=15$. Since $15$ is a quadratic residue mod $p$, only the cases
$p\equiv 1,2,4,8$ mod 15 are possible. Here the primitive characters $\psi_3=\LS -3$, $\psi_5$ as above, and $\psi_{15}=\LS -{15}$ occur.
Let $\EPS\in\{\pm 1\}$ and $p>5$ (which means $p\NDIV b$). Observe that $B_{\chi\psi_3}=-h(-3p)$ and $B_{\chi\psi_{15}}=-h(-15p)$.
We use the formulas (\ref{1.4}) again and obtain the following corollary.

\begin{corollary} % Corollary 2 %%%%%%%%%%%%%%%%%%%%%%%%%%%%%%%%%%%%%%%%%%%%%%%%%%%%%%%%%%%%%
\label{c2}
In the above setting,
\BD
  \sigma^2=T(15,p)+\frac{20h(-3p)-12\EPS\RE(B_{\chi\psi_5}(3-i))+15h(-15p)}{2(p-1)},
\ED
where $\EPS=1$ if $p\equiv 1$ mod $15$ and $\EPS=-1$ if $p\equiv 4$ mod $15$. Moreover,
\BD
  \sigma^2=T(15,p)+\frac{10h(-3p)+12\EPS\RE(B_{\chi\psi_5}(1-2i))+15h(-15p)}{2(p-1)},
\ED
where $\EPS=1$ if $p\equiv 2$ mod $15$ and $\EPS=-1$ if $p\equiv 8$ mod $15$.

\end{corollary} %%%%%%%%%%%%%%%%%%%%%%%%%%%%%%%%%%%%%%%%%%%%%%%%%%%%%%%%%%%%%%%%%%%%%%%%%%%%%

For instance, in the case $p=233$ ($\equiv 8$ mod $15$) we have $T(15,p)=539/29$. Since $h(-3p)=10$, $B_{\chi\psi_5}=-2-8i$, and $h(-15p)=36$, the second case of the corollary (with $\EPS=-1$) gives
$\sigma^2=539/29+107/58=1185/58$.

We expect that the double sum on the right-hand side of (\ref{1.3}) may result in arbitrarily complicated explicit formulas for a specific $b$ with many prime factors.

At this point one should say something about the natural density of the primes $p$ in question for a given $b$. This density is known only if one presupposes the Generalized Riemann Hypothesis.
Let $b$ be square-free, $b\equiv 2$ mod $4$ (so $b=10$ falls under this case).
Then the primes  $p$ such that $b$ has the order $q=p-1$ have the density
$A$, where $A$ is Artin's constant
\BD
  A=\prod_{l}\left(1-\frac 1{l(l-1)}\right)=0.3739558\LD,
\ED
$l$ running through the set of primes.
In our case $q=(p-1)/2$ the density is $3A/4$; see \cite[Th. 1]{Mo}.
In other words, one expects that more than 65\% of all primes below a large bound fall under the cases $q=p-1$ or $q=(p-1)/2$ for this type of $b$.

The connection between the digits of rational numbers and class number factors occurs in a number of papers; see \cite{Gi0}, \cite{Hi}, \cite{MuTh},
\cite{ChKr},  \cite{Mi}, \cite{PuSa}, and \cite{Sh}.

%%%%%%%%%%%%%%%%%%%%%%%%%%%%%%%%%%%%%%%%%%%%%
\section*{2. Proofs}
%%%%%%%%%%%%%%%%%%%%%%%%%%%%%%%%%%%%%%%%%%%%%

Let $q=(p-1)/2$ be the order of $b$ mod $p$. Hence $(c_1,\LD,c_q)$ is a period of the digit expansion (\ref{1.1}) of $1/p$. For an integer $k$ we define $(k)_p$ as the integer $j\in\{0,1,\LD,p-1\}$ that satisfies
$j\equiv k$ mod $p$. We put
\BD
  \theta_b(k)=(b(k)_p-(bk)_p)/p
\ED
for $k\in\Z$. Then the digit $c_j$ of $1/p$ has the form
\BD
  c_j=\theta_b(b^{j-1}),
\ED
for all $j\ge 1$; see \cite{Gi1}.
Put
\BD
  S=\sum_{j=1}^q c_j^2.
\ED
Let $m$ be the mean-value of the digits; see (\ref{1.1.1}).
We will use the identity
\BE
\label{2.1}
  \sigma^2=\frac{1}{q}S-m^2.
\EE
Since $b$ generates the set of quadratic residues mod $p$, the numbers $(b^{j-1})_p$, $j=1\LD,q$, run through the set $Q$ of quadratic residues mod $p$ in $\{1,\LD,p-1\}$, so
\BD
  S=\sum_{r\in Q}\theta_b(r)^2.
\ED

\MN
{\em Proof of Theorem \ref{t1}.} Let $p\equiv 3$ mod $4$, $p>3$. Since $-1$ is a quadratic nonresidue mod $p$, the numbers $(-b^{j-1})_p$, $j=1,\LD,q$, run through the set $N$ of quadratic nonresidues mod $p$
in $\{1,\LD,p-1\}$. Put
\BD
  c'_j=\theta_b(-b^{j-1}),
\ED
$j=1,\LD,q$, and
\BD
  S'=\sum_{j=1}^q c_j'^{\,2}\left(=\sum_{r\in N}\theta_b(r)^2\right).
\ED
Since $Q\cup N=\{1,\LD,p-1\}$,
\BD
  S+S'=\sum_{r=1}^{p-1}\theta_b(r)^2.
\ED
Therefore, formula (13) of \cite{Gi2} can be applied. It says
\BE
\label{2.2}
 S+S'=2bs(p,b)+\frac{(b-1)(2bp-3b-p+3)}6.
\EE
We use
\BD
 c_j'=b-1-c_j,\enspace j=1,\LD,q.
\ED
Accordingly,
\BD
 S'=S-2(b-1)\sum_{j=1}^q c_j+(b-1)^2(p-1)/2.
\ED
By \cite[Satz 11]{Gi1},
\BE
\label{2.4}
 \sum_{j=1}^q c_j=(b-1)(p-1)/4-(b-1)h(-p)/2.
\EE
We obtain
\BD
  S'=S+(b-1)^2h(-p).
\ED
Together with (\ref{2.2}), this gives
\BD
  2S=2bs(p,b)+\frac{(b-1)(2bp-3b-p+3)}6-(b-1)^2h(-p).
\ED
By (\ref{2.4}), the mean value $m$ of the digits equals
\BD
  m=\frac{b-1}2-\frac{b-1}{p-1}h(-p).
\ED
Since $q=(p-1)/2$,
we obtain from (\ref{2.1})
\BD
 \sigma^2=\frac{2bs(p,b)}{p-1}+\frac{(b-1)(2bp-3b-p+3)}{6(p-1)}-\frac{(b-1)^2}4-\frac{(b-1)^2}{(p-1)^2}h(-p)^2
\ED
and, thus, the desired result.
\STOP

\MN
{\em Proof of Theorem \ref{t2}.} Since $p\equiv 1$ mod 4, we have
\BD
  -1\equiv b^{(p-1)/4} \MB{ mod } p
\ED
and, therefore,
\BD
  c_{j+(p-1)/4}=b-1-c_j,\enspace j=1,\LD,(p-1)/4.
\ED
This implies (recall $q=(p-1)/2$)
\BD
  \sum_{j=1}^q c_j=(b-1)(p-1)/4.
\ED
Accordingly, the mean value of $c_1,\LD,c_q$ is
\BE
\label{2.6}
   m=(b-1)/2.
\EE

Let $l\in N$. Then
\BD
\frac lp=\sum_{j=1}^{\infty}c_j'p^{-j}\: \MB { with } c_j'=\theta_b(lb^{j-1});
\ED
see \cite{Gi2}. Put
\BD
  S'=\sum_{j=1}^q c_j'^{\,2}.
\ED
Since the numbers $(lb^{j-1})_p$ run through $N$,
we have (\ref{2.2}) with {\em this} definition of $S'$.
So $S+S'$ is made explicit in our sense --- but we also need an \BQ explicit\EQ expression for $S-S'$.

It is not hard to check that, for $k\in\{0,\LD,b-1\}$,
\BD
   \theta_b(r)=k\: \MB{ if, and only if }\: r\in \left(\frac{kp}b,\frac{(k+1)p}b\right).
\ED
Let $j$ run through $1,\LD,q$. Then
\BD
  |\{j;c_j=k\}|=|\{r\in Q;\theta_b(r)=k\}|=\left|Q\cap \left(\frac{kp}b,\frac{(k+1)p}b\right)\right|.
\ED
In the same way,
\BD
  |\{j;c_j'=k\}|=\left|N\cap \left(\frac{kp}b,\frac{(k+1)p}b\right)\right|.
\ED
Let
\BD
n_{b,k}=\left|\Z\cap\left(\frac{kp}b,\frac{(k+1)p}b\right)\right|
\ED
and
\BD
\delta_{b,k}=\left|Q\cap \left(\frac{kp}b,\frac{(k+1)p}b\right)\right|-\left|N\cap \left(\frac{kp}b,\frac{(k+1)p}b\right)\right|,
\ED
$k=0,\LD,b-1$. Since
\BD
 \left| Q\cap \left(\frac{kp}b,\frac{(k+1)p}b\right)\right|+ \left|N\cap \left(\frac{kp}b,\frac{(k+1)p}b\right)\right|=n_{b,k},
\ED
we obtain
\BD
  |\{j;c_j=k\}|=(n_{b,k}+\delta_{b,k})/2\: \MB{ and }\:|\{j;c_j'=k\}|=(n_{b,k}-\delta_{b,k})/2,
\ED
$j$ running through $1,\LD,q$. This implies
\BD
  S-S'=\sum_{k=0}^{b-1} k^2\cdot|\{j;c_j=k\}|-\sum_{k=0}^{b-1} k^2\cdot|\{j;c_j'=k\}|=\sum_{k=0}^{b-1}\delta_{b,k}k^2.
\ED
Combined with (\ref{2.2}), this identity gives
\BE
\label{2.6.0}
  2S=2bs(p,b)+\frac{(b-1)(2bp-3b-p+3)}6+\sum_{k=0}^{b-1}\delta_{b,k}k^2.
\EE

We have to deal with $\sum_{k=0}^q\delta_{b,k}k^2$ now. To this end we use the character $\chi=\LS -p$ again and define, for $k\in\{0,\LD,b\}$,
\BD
  \gamma_{b,k}=\sum_{1\le r\le kp/b}\chi(r).
\ED
Note that $\chi$ is an even character and $\gamma_{b,0}=0$. Then
\BD
  \delta_{b,k}=\gamma_{b,k+1}-\gamma_{b,k},\: k=0,\LD,b-1.
\ED
We obtain
\BD
  \sum_{k=0}^{b-1}\delta_{b,k}k^2=\sum_{k=0}^{b-1}(\gamma_{b,k+1}-\gamma_{b,k})k^2=\sum_{k=1}^{b}\gamma_{b,k}(k-1)^2-\sum_{k=1}^{b-1}\gamma_{b,k}k^2.
\ED
Since $\gamma_{b,b}=0$, we have
\BE
\label{2.6.1}
 \sum_{k=0}^{b-1}\delta_{b,k}k^2=-2\sum_{k=1}^{b-1}\gamma_{b,k}k+\sum_{k=1}^{b-1}\gamma_{b,k}.
\EE
Suppose that $d=(k,b)$. Then
\BD
  \gamma_{b,k}=\gamma_{b/d,k/d}.
\ED
For the sake of brevity, we write $d'$ instead of $b/d$. Then
\BE
\label{2.7}
  \sum_{k=1}^{b-1}\gamma_{b,k}k=\sum_{d\DIV b} d \sum_{\btop{l=1}{(l,d')=1}}^{d'-1}\gamma_{d',l}l.
\EE
Since $(l,d')=1$ we have, by formula (6) of \cite{Sz},
\BE
\label{2.8}
  \gamma_{d',l}=\frac{\chi(d')}{\phi(d')}\sum_{\psi\in X_{d'}^-}\OV{\psi}(-lp)B_{\chi\psi}.
\EE
Now formulas (\ref{2.7}) and (\ref{2.8}) give
\BD
\sum_{k=1}^{b-1}\gamma_{b,k}k=-\sum_{d\DIV b} d\frac{\chi(d')}{\phi(d')}\sum_{\psi\in X_{d'}^-}\OV{\psi}(p)B_{\chi\psi}\sum_{\btop{l=1}{(l,d')=1}}^{d'-1}\OV{\psi}(l)l.
\ED
Since
\BD
  \sum_{\btop{l=1}{(l,d')=1}}^{d'-1}\OV{\psi}(l)l=d'B_{\OV{\psi}},
\ED
we finally obtain
\BE
\label{2.10}
 \sum_{k=1}^{b-1}\gamma_{b,k}k=-b\sum_{d\DIV b}\frac{\chi(d)}{\phi(d)}\sum_{\psi\in X_{d}^-}\OV{\psi}(p)B_{\chi\psi}B_{\OV{\psi}}.
\EE
The sum
\BD
 \sum_{k=1}^{b-1}\gamma_{b,k}
\ED
is treated in the same way. We have
\BD
 \sum_{k=1}^{b-1}\gamma_{b,k}=-\sum_{d\DIV b}\frac{\chi(d')}{\phi(d')}\sum_{\psi\in X_{d'}^-}\OV{\psi}(p)B_{\chi\psi}\sum_{\btop{l=1}{(l,d')=1}}^{d'-1}\OV{\psi}(l).
\ED
But each $\psi\in X_{d'}^-$ is a non-trivial character, which means that the last sum on the right-hand side vanishes. Therefore,
\BE
\label{2.12}
  \sum_{k=1}^{b-1}\gamma_{b,k}=0.
\EE
In view of (\ref{2.1}), the theorem follows if we put (\ref{2.6}), (\ref{2.6.0}), (\ref{2.6.1}), (\ref{2.10}), and (\ref{2.12}) together.

\MN
This preprint has not undergone
peer review  or any post-submission improvements or corrections. The Version of Record of this article is
published in {\em The Ramanujan Journal}, and is available online at https://doi.org/10.1007/s11139-026-01406-5.

%%%%%%%%%%%%%%%%%%%%%%%%%%%%%%%%%%%%%%%%%%%%%%%%%%%%%
%%%%%%%%%%%%%%%%%%%%%%%%%%%%%%%%%%%%%%%%%%%%%%%%%%%%%%%%%%%%%%%%%%%%%%%%%%

\end{document}